\documentclass[12pt]{amsart}

\usepackage{latexsym}
\usepackage{amsmath}
\usepackage{amssymb}
\usepackage{amscd}
\usepackage{multicol}

\newtheorem{theorem}{\bf Theorem}[section]
\newtheorem{corollary}[theorem]{\bf Corollary}
\newtheorem{proposition}[theorem]{\bf Proposition}
\newtheorem{lemma}[theorem]{\bf Lemma}

\newtheorem{definition-theorem}[theorem]{\bf Theorem-Definition}

\def\bR{\mathbb{R}}
\def\bC{\mathbb{C}}
\def\bZ{\mathbb{Z}}

\def\D{{\mathcal D}}
\def\t{\mathfrak{t}}

\def\t{\frak{t}}

\def\h{\hbar}

\def\I{{\mathcal I}}
\def\E{{\mathcal E}}
\title[The  quantum cohomology ring of $G/B$]{A characterization of the quantum cohomology ring of $G/B$ and applications}

\date{\today}
\begin{document}
\begin{abstract} We show that the small quantum product of the
generalized flag manifold $G/B$ is a product operation on
$H^*(G/B)\otimes \bR[q_1,\ldots, q_l]$  uniquely determined by the
fact that it is a deformation of the cup product on $H^*(G/B)$, it
is commutative, associative, graded with respect to $\deg(q_i)=4$,
it satisfies a certain  relation (of degree two), and the
corresponding Dubrovin connection is flat. We deduce that it is
again the flatness of the Dubrovin connection which characterizes
essentially the solutions of the ``quantum Giambelli problem" for
$G/B$. 
This result gives new proofs of 
 the quantum Chevalley formula (see D.
Peterson [Pe] and Fulton and Woodward [Fu-Wo]), and of  Fomin, Gelfand and Postnikov's
description of the quantization map for $Fl_n$  (see [Fo-Ge-Po]).

\vspace{0.5cm}


\end{abstract}

\author[A.-L. Mare]{Augustin-Liviu  Mare}
\address{
Department of Mathematics and Statistics\\ University of
Regina \\ College West 307.14 \\ Regina SK, Canada S4S
0A2}
\email{mareal@math.uregina.ca}

\maketitle

\section{Introduction}

Let us consider the complex flag manifold $G/B$, where $G$ is a
connected, simply connected,  simple, complex Lie group and
$B\subset G$ a Borel subgroup. Let  $\t$ be the Lie algebra of a
maximal torus of a compact real form of $G$ and  $\Phi \subset
\t^*$ the corresponding set of roots.  Consider an arbitrary
 $W$-invariant inner product $\langle \ , \  \rangle$ on $\t$.
 To any root $\alpha$ corresponds the coroot
$$\alpha^{\vee}:= \frac{2\alpha}{\langle \alpha,
\alpha\rangle}$$ which  is an element of $\t$, by using the
identification of $\t$ and $\t^*$ induced by $\langle \ , \
\rangle$. If $\{\alpha_1, \ldots ,\alpha_l\}$ is a system of
simple roots then $\{\alpha_1^{\vee},\ldots, \alpha_l^{\vee}\}$ is
a  system of simple coroots.  Consider  $\{\lambda_1 ,\ldots ,
\lambda_l\} \subset \t^*$ the corresponding system of fundamental
weights, which are defined by
$\lambda_i(\alpha_j^{\vee})=\delta_{ij}$. The Weyl group $W$ is
the subgroup of $O(\t, \langle \ , \ \rangle )$ generated by the
reflections about the hyperplanes $\ker \alpha $,
 $\alpha\in \Phi^+$. It can be shown that
$W$ is in fact generated by the {\it simple reflections}
$s_1=s_{\alpha_1},\ldots, s_l=s_{\alpha_l}$ about the hyperplanes
$\ker\alpha_1, \ldots, \ker\alpha_l$. The {\it length} $l(w)$ of
$w$ is the minimal number of factors in a decomposition of $w$ as
a product of simple reflections. We denote by $w_0$ the longest
element of $W$.

Let $B^-\subset G$ denote the Borel subgroup opposite to $B$. To
each $w\in W$ we assign the {\it Schubert variety}
$X_w=\overline{B^-.w}$. The Poincar\'e dual of $[X_w]$ is an
element of $H^{2l(w)}(G/B)$, which is called the {\it Schubert
class}. The set $\{\sigma_w~|~w\in W\}$ is a basis of
$H^*(G/B)=H^*(G/B,\bR)$, hence  $\{\sigma_{s_1},\ldots,
\sigma_{s_l}\}$ is a basis of $H^2(G/B)$. A theorem of Borel [Bo]
says that the map
\begin{equation}\label{borel}
H^*(G/B)\to S(\t^*)/S(\t^*)^W=\bR[\{\lambda_i\}]/I_W\end{equation}
described by $ \sigma_{s_i} \mapsto [\lambda_i]$, $1\le i\le l$,
is a ring isomorphism (we are denoting by $S(\t^*)^W=I_W$ the
ideal of $S(\t^*)=\bR[\{\lambda_i\}]$ generated by the non-constant 
$W$-invariant polynomials). We will frequently identify $H^*(G/B)$
with the quotient ring from above.

To any $l$-tuple $d=(d_1,\ldots ,d_l)$ with $d_i\in \bZ$, $d_i\ge
0$ corresponds a {\it Gromov-Witten invariant}. This assigns to
any three Schubert classes $\sigma_u,\sigma_v, \sigma_w$ the
number denoted by $\langle \sigma_u|\sigma_v|\sigma_{w}\rangle_d$,
which
 counts the  holomorphic curves $\varphi :\bC P^1 \to
G/B$ such that $\varphi_*([\bC P^1])= d$ in $H_2(G/B)$ and
$\varphi(0)$, $\varphi(1)$ and $\varphi(\infty)$ are in general
translates of the Schubert varieties dual to $\sigma_u$,
$\sigma_v$, respectively $\sigma_{w}$. Let us consider the
variables $q_1,\ldots, q_l$. The {\it quantum cohomology ring} of
$G/B$ is the space $H^*(G/B)\otimes \bR[\{q_i\}]$ equipped with
the product $\circ$ which is  $\bR[\{q_i\}]$-linear and for any
two Schubert classes $\sigma_u, \sigma_v$, $u,v\in W$ we have
$$ \sigma_u \circ \sigma_v =\sum_{d=(d_1,\ldots ,d_l)\geq 0}
q^d\sum_{w\in W} (\sigma_u \circ \sigma_v)_d \sigma_w ,$$ $u,v\in
W$. Here $q^d$ denotes $q_1^{d_1}\ldots q_l^{d_l}$ and the
cohomology class $(\sigma_u \circ \sigma_v)_d$ is determined by
\begin{equation}\label{gw}\langle (\sigma_u \circ \sigma_v)_d,\sigma_w\rangle =
\langle \sigma_u|\sigma_v|\sigma_{w}\rangle_d,\end{equation} for
any $w\in W$. It turns out that the product $\circ$ is
commutative, associative and it is a deformation of the cup
product (by which mean that if we formally set $q_1=\ldots
=q_l=0$, then $\circ$ becomes the same as the cup product). If we
assign
$$\deg q_i=4,\quad 1\le i \le l,$$
then we also have the grading condition
$$\deg (a\circ b)=\deg a +\deg b,$$ for any two homogeneous elements $a, b$ of
$ H^*(G/B)\otimes
\bR[\{q_i\}] $. For more details about quantum cohomology we refer
the reader to Fulton and Pandharipande [Fu-Pa].

The {\it Dubrovin connection} attached to the quantum product
defined above is a connection\footnote{More precisely, a family of
connections depending on the parameter $\h\in \bR\setminus
\{0\}$.} $\nabla ^{\hbar}$ on the trivial vector bundle
$H^*(G/B)\times H^2(G/B)\to H^2(G/B)$ defined as follows: Denote
by $t_1,\ldots, t_l$ the coordinates on $H^2(G/B)$ induced by the
basis $\sigma_{s_1}, \ldots, \sigma_{s_l}$. Consider the 1-form
$\omega$ on $H^2(G/B)$ with values in ${\rm End}(H^*(G/B))$ given
by
$$\omega_t(X,Y) = X\circ Y,$$ for $t=(t_1,\ldots, t_l) \in
H^2(G/B)$, $X\in H^2(G/B)$ and $Y\in H^*(G/B)$, where the
convention $$q_i=e^{t_i}, \quad 1\le i\le l$$ is in force. Finally
set
$$\nabla^\h  = d +\frac{1}{\h} \omega.$$
Note that the 1-form $\omega$  can be expressed as
$$\omega = \sum_{i=1}^l \omega_i dt_i,$$
where $\omega_i$ denotes the matrix of the operator
$\sigma_{s_i}\circ$ on $H^*(G/B)$ with respect to the basis
consisting of the Schubert classes. The following result is
well-known (cf. [Du]):

\begin{lemma} The Dubrovin connection $\nabla^\h$ is flat for any
$\h\in \bR\setminus \{0\}$, i.e. we have
\begin{equation}\label{zero}d\omega = \omega\wedge \omega =
0\end{equation}
\end{lemma}
\begin{proof} The fact that $d\omega=0$ follows from
$$\frac{\partial}{\partial t_i} \omega_j =\frac{\partial}{\partial t_j} \omega_i
,$$ which is equivalent to
$$d_i(\sigma_{s_j}\circ\sigma_w)_d = d_j
(\sigma_{s_i}\circ\sigma_w)_d$$ for any $w\in W$ and any
$d=(d_1,\ldots, d_l)$, hence, by (\ref{gw}), to
$$d_i\langle \sigma_{s_j}|\sigma_w |\sigma_v\rangle_d=
d_j\langle \sigma_{s_i}|\sigma_w |\sigma_v\rangle_d.$$ The latter
equality follows from the ``divisor property" (see [Fu-Pa,
equation (40)] for a more general version of this formula):
$$\langle \sigma_{s_j}|\sigma_w |\sigma_v\rangle_d =
d_j \langle \sigma_w |\sigma_v\rangle_d.$$ The  equality
$\omega\wedge \omega=0$, i.e. $\omega_i\omega_j =
\omega_j\omega_i$,  $1\le i,j\le l$, follows from the fact that
the product $\circ$ is commutative and associative.
\end{proof}

Another important property of the quantum product which is of
interest for us  is that we have the relation:
\begin{equation}\label{relation}\sum_{i,j=1}^l\langle \alpha_i^{\vee},
\alpha_j^{\vee}\rangle \sigma_{s_i}\circ\sigma_{s_j} =
\sum_{i=1}^l\langle \alpha_i^{\vee}, \alpha_i^{\vee}\rangle
q_i.\end{equation} In order to prove this we  take into account
that: \begin{itemize} \item we have (see [Kim] or [Ma1, Lemma
3.2])
$$\sigma_{s_i}\circ \sigma_{s_j} = \sigma_{s_i}\sigma_{s_j}
+\delta_{ij}q_j$$  \item the polynomial $\sum_{i,j=1}^l\langle
\alpha_i^{\vee}, \alpha_j^{\vee}\rangle \lambda_{i}\lambda_{j}\in
S(\t^*)$ is $W$-invariant (being just the squared norm on $\t$);
hence, according to (\ref{borel}), the following relation holds in
$H^*(G/B)$:
 $$\sum_{i,j=1}^l\langle \alpha_i^{\vee},
\alpha_j^{\vee}\rangle \sigma_{s_i}\sigma_{s_j}=0.$$
\end{itemize}
The goal of this paper is to show that the quantum product for
$G/B$ is essentially determined by the equations (\ref{zero}) (the
flatness of the Dubrovin connection) and (\ref{relation}) (the
degree two relation). More precisely, we will prove that:

\begin{theorem}\label{main} Let $\star$ be a product on the space $H^*(G/B)\otimes
\bR[\{q_i\}]$ which is commutative, associative,  is a deformation
of the cup product (in the sense defined above),   satisfies the
condition $\deg(a\star b)=\deg a +\deg b$, for $a,b$ homogeneous
elements of $H^*(G/B)\otimes \bR[\{q_i\}]$, with respect to the
grading $\deg q_i=4$, and

\begin{itemize}
\item[(a)] the Dubrovin connection $\nabla ^\h = d+\frac{1}{\h}
\omega$, with $\omega(X,Y)= X\star Y$ is flat. In other words, if $\omega_k$ is
the matrix of the $ \bR[\{q_i\}]$-linear endomorphism
$\sigma_{s_k}\star$ of $H^*(G/B)\otimes \bR[\{q_i\}]$ with respect to the Schubert basis,
then we have
$$\frac{\partial}{\partial t_i} \omega_j =\frac{\partial}{\partial t_j} \omega_i
$$ for all $1\le i,j\le l$ (the convention $q_i=e^{t_i}$ is in force). \item[(b)] we have
$$\sum_{i,j=1}^l\langle \alpha_i^{\vee}, \alpha_j^{\vee}\rangle
\sigma_{s_i}\star\sigma_{s_j} = \sum_{i=1}^l\langle
\alpha_i^{\vee}, \alpha_i^{\vee}\rangle q_i.$$
\end{itemize}
Then $\star$ is the quantum product $\circ$.
\end{theorem}

The proof will be done in section 2. There are two corollaries we
would like to deduce from this theorem. The first one is a
characterization of the quantum Giambelli polynomials in terms of
the flatness of the Dubrovin connection. More precisely, let us
denote by $QH^*(G/B)$ the quotient ring $\bR[\{\lambda_i\},
\{q_i\}]/\langle R_1,\ldots, R_l \rangle$, where $R_1, \ldots,
R_l$ are the quantum deformations in the quantum cohomology ring
$(H^*(G/B)\otimes \bR[\{q_i\}], \circ)$ of the fundamental
homogeneous generators of $S(\t^*)^W$  ($R_1,\ldots, R_l$ have
been determined explicitly  by B. Kim in [Kim];  we will present
in section 2 a few more details about that). For any $c\in
\bR[\{\lambda_i\},\{q_i\}]$ we denote by $[c]_q$ the coset of $c$
in $QH^*(G/B)$. The map $\sigma_{s_i}\mapsto [\lambda_i]_q$
induces a tautological isomorphism
\begin{equation}\label{tautological}(H^*(G/B)\otimes
\bR[\{q_i\}],\circ)\simeq QH^*(G/B).\end{equation} Finding for
each $w\in W$ a polynomial $\hat{c}_w\in
\bR[\{\lambda_i\},\{q_i\}]$ whose coset in $QH^*(G/B)$ is the
image of $\sigma_w$ --- in other words, solving the quantum
Giambelli problem --- would lead to a complete knowledge of the
quantum cohomology of $G/B$. We are looking for conditions which
determine the polynomials $\hat{c}_w$. First of all, let us
consider for each $w\in W$ a polynomial\footnote{These are
 solutions of the classical Giambelli problem for $G/B$. Such polynomials have
been constructed for instance by Bernstein, I. M. Gelfand and S.
I. Gelfand in [Be-Ge-Ge].} $c_w\in \bR[\{\lambda_i\}]$  whose
coset corresponds to $\sigma_w$ via the isomorphism (\ref{borel}).
There are two natural conditions that we impose to the polynomials
$\hat{c}_w$:
\begin{equation}\label{deg}\deg \hat{c}_w =\deg c_w\end{equation}
with respect to the grading $\deg \lambda_i=2$, $\deg q_i=4$, and
\begin{equation}\label{hat}\hat{c}_w|_{({\rm all} \ q_i \ =0)}=c_w.\end{equation}
Whenever the conditions (\ref{deg}) and (\ref{hat}) are satisified,
the cosets $[\hat{c}_w]_q$, $w\in W$, are a basis of $QH^*(G/B)$ over
$\bR[\{q_i\}]$. Consider the 1-form $$\omega=\sum_{i=1}^l\omega_i dt_i,$$
where $\omega_i$ is the matrix of multiplication of $QH^*(G/B)$ by $[\lambda_i]_q$
with respect to the latter basis. We can prove that:
\begin{corollary}\label{proper}
Let $\hat{c}_w$, $w\in W$, be polynomials in $\bR[\{\lambda_i\}, \{q_i\}]$
which satisfy the properties
(\ref{deg}) and (\ref{hat}).  Then the image of $\sigma_w$ by the isomorphism
(\ref{tautological}) is $[\hat{c}_w]_q$ for all $w\in W$ if and only if
the connection $$\nabla^{\hbar} =d+\frac{1}{\hbar}\omega$$ is flat
for all $\hbar\in \bR\setminus\{0\}$. The latter condition reads
$$\frac{\partial}{\partial t_i} \omega_j =\frac{\partial}{\partial t_j} \omega_i,$$
for all $1\le i,j\le l$.
\end{corollary}
\begin{proof} Consider the $\bR[\{q_i\}]$-linear isomorphism\footnote{This is what
Amarzaya and Guest [Am-Gu] call a ``quantum evaluation map".}
$$\delta : QH^*(G/B)\to H^*(G/B)\otimes \bR[\{q_i\}]=\bR[\{\lambda_i\}, \{q_i\}]/(I_W\otimes
\bR[ \{q_i\}])$$  determined by
\begin{equation}\label{delta}\delta[\hat{c}_w]_q = [c_w],\end{equation} for all $w\in W$.
Define
 the product $\star$ on
$H^*(G/B)\otimes \bR[\{q_i\}]$ by
$$x\star y = \delta(\delta^{-1}(x) \delta^{-1}(y)),$$
$x,y \in H^*(G/B)\otimes \bR[\{q_i\}]$. The product is
commutative, associative, it is a deformation of the cup product
on $H^*(G/B)$, and it satisfies  $\deg(a\star b)=\deg a +\deg b$,
where $a,b\in H^*(G/B)\otimes \bR[\{q_i\}]$ are homogeneous
elements. The map $\delta$ is obviously a ring isomorphism between
$QH^*(G/B)$ and $(H^*(G/B)\otimes \bR[\{q_i\}],\star)$. In
particular, the following degree  two relation holds:
$$\sum_{i,j=1}^l\langle \alpha_i^{\vee}, \alpha_j^{\vee}\rangle
[\lambda_i]\star[\lambda_j] = \sum_{i=1}^l\langle \alpha_i^{\vee},
\alpha_i^{\vee}\rangle q_i.$$ Moreover, the matrix of
$[\lambda_i]\star$ on $H^*(G/B)\otimes \bR[q_1,\ldots, q_l]$ with
respect to the Schubert basis $\{[c_w]: w\in W\}$ is just
$\omega_i$. So if the connection $\nabla^{\hbar}$ is flat for all
$\hbar$, then, by Theorem \ref{main}, the products $\star$ and
$\circ$ are the same. This implies that $\delta$ is just the
isomorphism (\ref{tautological}). The conclusion follows from the
definition (\ref{delta}) of $\delta$.
\end{proof}

Corollary \ref{proper} will be used in section \ref{last} in order to recover the
``quantization via standard monomials" theorem of Fomin, Gelfand, and Postnikov
(see [Fo-Ge-Po, Theorem 1.1]).   

Our second  application of Theorem \ref{main} concerns 
the combinatorial quantum product
$\star$ on $H^*(G/B)\otimes \bR[\{q_i\}]$, which has been constructed in [Ma4]. 
By definition, this
product satisfies the following {\it quantum Chevalley formula}:
$$\sigma_{s_i}\star \sigma_w = \sigma_{s_i}\sigma_w + \sum_{}
\lambda_i(\alpha^{\vee})\sigma_{ws_{\alpha}},$$ for $1\le i\le l$,
$w\in W$. Here the sum runs over all positive roots $\alpha$ with
the property that $l(ws_{\alpha})=l(w)-2{\rm
height}(\alpha^{\vee}) +1$, where we consider the expansion
$\alpha^{\vee}=m_1\alpha_1^{\vee}+\ldots + m_l\alpha_l^{\vee}$,
$m_j\in \bZ$, $m_j\ge 0$ and denote
$${\rm
height}(\alpha^{\vee})=m_1+\ldots +m_l ,\quad
q^{\alpha^{\vee}}=q_1^{m_1}\ldots q_l^{m_l}.$$ We have also showed
in [Ma4]  that $\star$ satisfies all hypotheses of Theorem
\ref{main}. We deduce:

\begin{corollary} The combinatorial and actual quantum products
coincide. Consequently, the quantum product $\circ$ satisfies the
 quantum Chevalley formula:
\begin{equation}\label{chevalley}\sigma_{s_i}\circ \sigma_w = \sigma_{s_i}\sigma_w +
\sum_{l(ws_{\alpha})=l(w)-2{\rm height}(\alpha^{\vee}) +1}
\lambda_i(\alpha^{\vee})\sigma_{ws_{\alpha}},\end{equation} for
$1\le i\le l$, $w\in W$.
\end{corollary}

{\bf Remarks.} 1. The formula (\ref{chevalley}) plays a crucial
role in the study of the quantum cohomology algebra of $G/B$, as
this  is generated over $\bR[q_1,\ldots, q_l]$ by the degree 2
Schubert classes $\sigma_{s_1},\ldots, \sigma_{s_l}$. The formula
was announced by D. Peterson in [Pe]. A rigorous
intersection-theoretic proof has been given by W. Fulton and C.
Woodward in [Fu-Wo]. It is one of the aims of our paper to give an
alternative, conceptually new, proof of this formula.

2. The tool we will be using in the proof of Theorem \ref{main} is
the notion of   $\D$-module, in the spirit of Guest [Gu], Amarzaya
and Guest [Am-Gu], and Iritani [Ir]. More precisely, we will show
that the $\D$-modules associated in Iritani's manner to the
products $\circ$ and $\star$ are isomorphic, and then we conclude
by using a certain uniqueness argument of Amarzaya and Guest
[Am-Gu] (for more details, see the next section).

{\bf Acknowledgements.} I would like to thank Jost Eschenburg and
Martin Guest for discussions on the topics contained in this paper.

\section{$\D$-modules and quantum cohomology}

The goal of this section is to give a proof of Theorem \ref{main}.

We denote by $\D$ the Heisenberg algebra, by which we mean the
associative $\bR[\h]$-algebra generated by $Q_1,\ldots, Q_l,$
$P_1,\ldots, P_l$, subject to  the relations
\begin{equation}\label{comm}[Q_i,Q_j]=[P_i,P_j]=0, \quad [P_i,Q_j] =\delta_{ij}
\h Q_j,\end{equation} $1\le i,j\le l$. It becomes a graded algebra
with respect to the assignments \begin{equation}\label{degree}\deg
Q_i=4, \quad \deg P_i=\deg \h=2.\end{equation} Note that any
element $D$ of $\D$ can be  written uniquely as an
$\bR[\h]$-linear combination of monomials of type $Q^IP^J$.

A concrete realization of $\D$ can be obtained by putting
$Q_i=e^{t_i}$ and $P_i=\h\frac{\partial}{\partial t_i}$, $1\le
i\le l$. We will be interested in certain elements of $\D$ which
arise in connection with the Hamiltonian system of Toda lattice
type corresponding to the coroots of $G$, namely the first quantum
integrals of motion of this system. Those are homogeneous elements
$D_k=D_k(\{Q_i\},\{P_i\},\h)$ of $\D$, $1\le k\le l$, which
commute with
$$D_1=\sum_{i,j=1}^l\langle \alpha_i^{\vee},
\alpha_j^{\vee}\rangle P_iP_j - \sum_{i=1}^l\langle
\alpha_i^{\vee}, \alpha_i^{\vee}\rangle Q_i$$ and also satisfy the
property that $D_k(\{0\}, \{\lambda_i\},0)$, $1\le k\le l$, are
just the fundamental homogeneous $W$-invariant polynomials (for
more details concerning the differential operators $D_1,\ldots,
D_l$ we address the reader to [Ma3]). We will denote by $\I$ the
left sided ideal of $\D$ generated by $D_1,\ldots, D_l$.

  Let $\star$ be  a product on $H^*(G/B)\otimes
\bR[\{q_i\}]$ which satisfies the hypotheses of Theorem
\ref{main}. Let us denote by $E$ the $\D$-module (i.e. vector
space with an action of the algebra $\D$) $H^*(G/B)\otimes
\bR[\{q_i\},\h]$ defined by
$$Q_i.a= q_i a,\quad P_i.a = \sigma_{s_i}\star a
+\h q_i\frac{\partial}{\partial q_i} a,$$ $1\le i\le l$, $a\in
H^*(G/B)\otimes \bR[\{q_i\},\h]$. The isomorphism type of the
$\D$-module $E$ corresponding to $\star$ is uniquely determined by
the hypotheses of Theorem \ref{main}, as the following proposition
shows:

\begin{proposition}\label{isomorphism} If $\star$ is a product with
the properties stated in Theorem \ref{main}, then the map
$\phi:\D\to H^*(G/B)\otimes \bR [\{q_i\}, \h]$ given by
 $$f(\{Q_i\},\{P_i\},\h) \stackrel{\phi}{\mapsto}  f(\{Q_i\}, \{P_i\}, \h).1=
f(\{q_i\},\{\sigma_{s_i}\star +\h q_i\frac{\partial}{\partial
q_i}\},\h).1$$ is surjective and induces an isomorphism of
$\D$-modules
\begin{equation}\label{isom}\D/\I \simeq E,\end{equation}
where $\I$ is the left sided ideal of $\D$ generated by the
quantum integrals of motion of the Toda lattice (see above).
\end{proposition}

\begin{proof}  We will use the grading on
$H^*(G/B)\otimes \bR[\{q_i\}, \h]$ induced by the usual grading on
$H^*(G/B)$, $\deg q_i=4$ and $\deg \h=2$. Combined with the
grading defined by (\ref{degree}), this makes $\phi$ into a degree
preserving map (more precisely, it maps a homogeneous element of
$\D$ to a homogeneous element of the same degree in
$H^*(G/B)\otimes \bR[\{q_i\}, \h]$).

 Let us prove first the surjectivity stated in our theorem.
It is sufficient to show that any homogeneous element $a\in
H^*(G/B)\otimes \bR[\{q_i\}, \h]$ can be written as $f(\{Q_i\},
\{P_i\}, \h).1$. We proceed by induction on $\deg a$. If $\deg
a=0$, everything is clear.  Now consider $a\in H^*(G/B)\otimes
\bR[\{q_i\}, \h]$ a homogeneous element of degree at least 2. By a
result of Siebert and Tian [Si-Ti], we can express
$$a=g(\{q_i\},\{\sigma_{s_i}\star\}, \h)$$ for a certain polynomial $g$. We
have
$$a-g(\{Q_i\},\{ P_i\}, \h).1=a-g(\{q_i\},\{\sigma_{s_i}\star +
\h q_i\frac{\partial}{\partial q_i}\} , \h).1 =\h b,$$ where $b\in
H^*(G/B)\otimes \bR[\{q_i\}, \h]$ is homogeneous of degree $\deg
a-2$ or it is zero. We use the induction hypothesis for $b$.

We proved in [Ma3] (see the proof of Lemma 4.5) that the
generators $D_k=D_k(\{Q_i\},\{P_i\},\h)$, $1\le k\le l$, of the
ideal  $\I$ satisfy
 \begin{equation}\label{homog} D_k(\{Q_i\},
\{P_i\},\h).1=0.
\end{equation} If we let $\h$ approach $0$ in
(\ref{homog}) we obtain the  relations
\begin{equation}\label{dek}D_k(\{q_i\},
\{\sigma_{s_i}\star\},0)=0,\end{equation} $1\le k\le l$. They
generate the whole ideal of relations in the ring
$(H^*(G/B)\otimes \bR[\{q_i\}],\star)$.

We need to show that if $D$ is an element of $\D$ with the
property that
\begin{equation}\label{de}D(\{Q_i\},\{P_i\},\h).1=0\end{equation} then
$D\in\I$. Because the map $\phi$ is degree preserving, we may
assume that $D$ is homogeneous and proceed by induction on $\deg
D$. If $\deg D=0$, i.e. $D$ is constant, then (\ref{de}) implies $D=0$, hence
$D\in \I$.  It now follows the
induction step. From
$$D.1=D(\{q_i\},\{\sigma_{s_i}\star + \h q_i\frac{\partial}{\partial q_i}\},
\h).1=0,$$ for all $\h$, we deduce the relation
$D(\{q_i\},\{\sigma_{s_i}\star\}, 0)=0$ in the ring
$(H^*(G/B)\otimes \bR[\{q_i\}],\star)$. Consequently we have the
following polynomial identity
$$D(\{q_i\},\{\lambda_i\},0) =
\sum_kf_k(\{q_i\},\{\lambda_i\})D_k( \{q_i\}, \{\lambda_i\}, 0),$$
for certain polynomials $f_k$. By using the commutation relations
(\ref{comm}), we obtain the following identity in $\D$:
\begin{align}D(\{Q_i\},\{P_i\},0) &\equiv \sum_k
f_k(\{Q_i\},\{P_i\})D_k(\{Q_i\},\{P_i\},0) \ {\rm mod} \ \h
\nonumber
\\ {} & \equiv \sum_k f_k(\{Q_i\},\{P_i\})D_k(\{Q_i\},\{P_i\},\h) \ {\rm
mod} \ \h.\nonumber \end{align} In other words,
$$D(\{Q_i\}, \{P_i\}, \h)= \sum_k f_k(\{Q_i\},\{P_i\})D_k(\{Q_i\},\{P_i\},\h)
+\h D'(\{Q_i\},\{P_i\}, \h),$$ for a certain $D'\in\D$, with $\deg
D'<\deg D$. From (\ref{dek}) and (\ref{de}) we deduce that
$$D'(\{Q_i\},\{P_i\},\h).1=0$$
Since $\deg D'<\deg D$, we only have to use the induction
hypothesis for $D'$ and get to the desired conclusion.

\end{proof}

Note that (\ref{isom}) is also an isomorphism of
$\bR[\{Q_i\},\h]$-modules. Since the actual quantum product
$\circ$ satisfies the hypotheses of Theorem \ref{main}, we deduce
that the dimension of $\D/\I$ as an $\bR[\{Q_i\},\h]$-module
equals $|W|$. Let us consider the ``standard monomial basis"
$\{[C_w]  : w\in W\}$ of $\D/\I$ over $\bR[\{Q_i\},\h]$ with
respect to a choice of a Gr\"obner basis of the ideal $\I$ (for
more details, see Guest [Gu, section 1] and the references
therein). Any $C_w$ is a monomial in $P_1,\ldots,P_l$ and the
cosets of the monomials
$$c_w=C_w(\lambda_1,\ldots,\lambda_l),\quad w\in W$$ in
$H^*(G/B)=S(\t^*)/S(\t^*)^W =\bR[\{\lambda_i\}]/I_W$ are a basis.
We will need the following result (our proof relies on an idea of
Amarzaya and Guest [Am-Gu]):
\begin{proposition}\label{unique} There exists a unique basis $\{[\bar{C}_w]:w\in W\}$ of
$\D/\I$ over $\bR[\{Q_i\},\h]$ with the following properties:
\begin{itemize}
\item[(i)] for all $w\in W$ the element $\bar{C}_w=\bar{C}_w(\{Q_i\},\{P_i\},\h)$ of $\D$ is
homogeneous of degree $2\deg c_w$  with respect to the grading
defined by (\ref{degree})
\item[{}]
\item[(ii)] for all $w\in W$ we have
$$\bar{C}_w(\{0\},\{\lambda_i\},\h) \equiv c_w\ {\rm mod}\ I_W;$$ in
particular $\bar{C}_w(\{0\},\{\lambda_i\},\h)  {\rm mod}\ I_W$ is
independent of $\hbar$
\item[{}]
\item[(iii)]   the elements
$(\bar{\Omega}^i_{vw})_{v,w\in W}^{1\le i\le l}$  of
$\bR[Q_1,\ldots, Q_l,\h]$ determined by
$$P_i[\bar{C}_w]=\sum_{v\in W}\bar{\Omega}^i_{vw}[\bar{C}_v],$$
 are independent of  $\h$. \end{itemize}
\end{proposition}
\begin{proof} In order to show that such a basis exists, we
consider the isomorphism $$\phi:\D/\I \to H^*(G/B)\otimes
\bR[\{q_i\},\h]$$ induced by the actual quantum product $\circ$
via Proposition \ref{isomorphism}. The basis $\{[c_w]:w\in W\}$ of
the right hand side induces the basis
$\{[\bar{C}_w]=\phi^{-1}([c_w]): w\in W\}$ of $\D/\I$ over
$\bR[\{Q_i\}, \h]$.  It is obvious that the latter basis satisfies
(i) and (iii).  In order to show that it also satisfies (ii), we
consider the following commutative diagram:
$$\D/\I \stackrel{\phi}{\longrightarrow} H^*(G/B)\otimes
\bR[\{q_i\}, \h]$$
$$\psi_1 \searrow \    \ \ \ \ \ \ \ \ \ \ \ \ \  \  \swarrow   \psi_2 \ \ \ \ \ \ \ \ {}$$
$$H^*(G/B) \otimes \bR[\h]\ \ \ \ \ \ \ \ \ {}$$
where $\psi_2$ is the canonical projection and $\psi_1:\D/\I \to
H^*(G/B) \otimes \bR[\h] =(\bR[\{\lambda_i\}]/I_W) \otimes
\bR[\h]$ is given by
$$[D(\{Q_i\},\{P_i\}, \h)] \mapsto [D(\{0\}, \{\lambda_i\}, \h)].$$
Note that $\psi_1$ is well defined, as for any $k=1,2,\ldots ,l$,
the polynomial  $D_k(\{0\}, \{\lambda_i\}, \h)$ is independent of
$\h$, being equal to $u_k$, the $k$-th fundamental  $W$-invariant
polynomial (see [Ma 2, section 3]). We observe that
$$[\bar{C}_w(\{0\}, \{\lambda_i\},
\h)]=\psi_1 [\bar{C}_w]=\psi_2[c_w]=[c_w],$$ hence condition (ii)
is satisfied.

In order to show that there exists at most one such basis, we will
use a construction of Amarzaya and Guest [Am-Gu]. Let
$\{[\bar{C}_w]: w\in W\}$ be a basis of $\D/\I$ with the
properties (i), (ii) and (iii). We can write
\begin{equation}\label{bar}\bar{C}_w \equiv \sum_{v\in W} U^{vw}C_v \ {\rm
mod} \ \I\end{equation} with $U^{vw}\in\bR[\{Q_i\},\h]$. Decompose
the matrix $U=(U^{vw})_{v,w\in W}$ as $$U=U_0+\h U_1 +\ldots +\h^k
U_k,$$ where
 $U_0,\ldots, U_k$ have entries in
$\bR[\{Q_i\}]$. Let us apply $\psi_1$ to both sides of equation
(\ref{bar}) and deduce that in $H^*(G/B)\otimes \bR[\h]
=(\bR[\{\lambda_i\}]/I_W) \otimes \bR[\h]$ we have that
$$[c_w] = \sum_{v\in W} U^{vw}|_{({\rm all} \ Q_i =0)}[c_v],$$
for all $w\in W$. This implies $$U^{vw}|_{({\rm all} \ Q_i
=0)}=\delta_{vw},$$ where $\delta_{vw}$ is the Kroenecker delta.
On the other hand, because any $\bar{C}_w$, $C_v$, $v,w\in W$, as
well as any generator $D_i$ of $\I$ is homogeneous, we deduce that
each $U^{vw}$ is homogeneous. We are  led  to the following
property of the matrices $U_j$:
\begin{itemize}
\item[(a)] besides the diagonal of $U_0$, which is $I$, the
entries of $U_0, U_1, \ldots, U_k$  are homogeneous polynomials
with no degree zero term in $Q_1,\ldots, Q_l$
\end{itemize}

 Let us choose an ordering of $W$ which is
increasing with respect to $\deg c_w$. In this way, the set
$\{[c_w]:w\in W\}$ is a  basis of $H^*(G/B)$ consisting of $s_0=1$
elements of degree $0$, followed by $s_1$ elements of degree $2$,
$\ldots$ , followed by $s_m$ elements of degree $2m=\dim G/B$. All
matrices involved here appear as block matrices of the type
$A=(A_{\alpha\beta})_{1\le\alpha,\beta \le m}$. We will say that a
block matrix $A=(A_{\alpha\beta})_{1\le\alpha,\beta \le m}$ is
$r$-{\it triangular} if $A_{\alpha\beta}=0$ for all $\alpha,\beta$
with $\beta -\alpha <r$. From the homogeneity of $U^{vw}$
mentioned above and the fact that $\deg Q_1=\ldots =\deg Q_l=4$,
we deduce:
\begin{itemize}
\item[(b)] the block matrix  $U_0-I$ is 2-triangular
\item[(c)] for
any $1\le j\le k$, the block matrix $U_j$ is $(j+2)$-triangular.
\end{itemize}
In particular we can assume that $k=m-2$, hence
\begin{equation}\label{u}U=U_0+\h U_1 +\ldots +\h^{m-2} U_{m-2},\end{equation}

Consider the matrix $\Omega^i=(\Omega^i_{vw})_{v,w\in W}$
determined by \begin{equation}\label{pi}P_i[C_w]=\sum_{v\in
W}\Omega^i_{vw}[C_v].\end{equation} As before, each
$\Omega^i_{vw}$ is an element of $\bR[Q_1,\ldots, Q_l, \h]$ which
is homogeneous with respect to the grading given by
(\ref{degree}). Also, if we apply $\psi_1$ on both sides of the
equation (\ref{pi}), we deduce that in $H^*(G/B)\otimes \bR[\h]
=(\bR[\{\lambda_i\}]/I_W) \otimes \bR[\h]$ we have
$$[\lambda_i][c_w]=\sum_v \Omega^i_{vw}|_{({\rm all} \ Q_i
=0)}[c_v].$$ This shows that \begin{equation}\label{omega}
\Omega^i_{vw}|_{({\rm all} \ Q_j =0)} \ {\rm is \ independent \ of
\ } \h, {\rm for  \ all} \ v,w\in W, 1\le i\le l\end{equation}

From here on, it will be more convenient to work with the
realization of $\D$ given by $Q_i=e^{t_i}$ and
$P_i=\h\frac{\partial}{\partial t_i}$, $1\le i\le l$. Then
$\Omega^i$ become matrices whose coefficients are homogeneous
polynomials in $e^{t_1},\ldots, e^{t_l}$, and $\hbar$. Let us
consider the 1-form
\begin{equation}\label{omega0}\Omega =\sum_{i=1}^l \Omega_i dt_i.\end{equation} We  decompose it as
$$\Omega= \omega +\h\theta^{(1)}+\ldots +
\h^p\theta^{(p)}.$$ From the homogeneity of the entries of
$\Omega_i$, as well as from (\ref{omega0}) we deduce that:
\begin{itemize}
\item[(d)] the block matrix $\omega$ is $(-1)$-triangular
\item[(e)] the block matrix $\theta^{(j)}$ is $(j+1)$-triangular, for any $1\le j \le
p$.
\end{itemize}
In particular we can assume that $p=m-1$, hence
\begin{equation}\label{omega}\Omega= \omega +\h\theta^{(1)}+\ldots +
\h^{m-2}\theta^{(m-2)}.\end{equation}  Now consider the matrix
$\bar{\Omega}^i=(\bar{\Omega}^i_{vw})_{v,w\in W}$ determined by
$$P_i[\bar{C}_w]=\sum_{v\in W}\bar{\Omega}^i_{vw}[\bar{C}_v].$$
Note that if $p\in \D$ is a polynomial $p(e^{t_1}, \ldots,
e^{t_l})$, then we have
$$\h\frac{\partial}{\partial t_i} \cdot p = p\cdot \h\frac{\partial}{\partial t_i}
+\h\frac{\partial}{\partial t_i}(p).$$ By using this, we can
easily deduce from (\ref{bar}) that
$$\bar{\Omega}^i= U^{-1}\Omega^i U+\h U^{-1}\frac{\partial}{\partial
t_i} U.$$ Thus the 1-form $\bar{\Omega} =\sum_{i=1}^l
\bar{\Omega}_i dt_i$ is given by
$$\bar{\Omega} = U^{-1}\Omega U +\h U^{-1}dU.$$
Condition (iii) reads $\bar{\Omega}$ is independent of $\h$.  From
(\ref{u}) and (\ref{omega}) we can see that this is equivalent to
$$U^{-1}\Omega U +\h U^{-1}dU =U_0^{-1}\omega U_0$$
and further to
\begin{equation}\label{diffeq}\Omega U +\h dU =UU_0^{-1}\omega U_0.\end{equation}
We will prove the following claim

{\it Claim.} For a given $\Omega$ of the type (\ref{omega}) with
the properties  (d) and (e), the system (\ref{diffeq})  has at
most one solution $U$ of the type (\ref{u}) with $U_j$ satisfying
(a) and (b).

It is obvious that the claim implies that there exists at most one
basis $\{[\bar{C}_w]:w\in W\}$ with the properties (i), (ii) and
(iii), and the proof is complete.

In order to prove the claim, let us write
$$U=(I +\h V_1 + \h^2 V_2 +\ldots +\h^{m-2}V_{m-2})V_0$$
where $V_0=U_0$, $V_1=U_1U_0^{-1},\ldots ,
V_{m-2}=U_{m-2}U_0^{-1}$. Note that (a), (b) and (c) from above
imply:
\begin{itemize}
\item[(f)]  $V_0$ is a block matrix whose diagonal is $I$,
such that $V_0-I$ is $2$-triangular,  and all entries of $V_0$
which are not on the diagonal are polynomials with no degree zero
term in $e^{t_1},\ldots, e^{t_l}$,
\item[(g)] $V_0^{-1}$ is an upper triangular matrix, its diagonal
is $I$, and all entries of $V_0^{-1}$ which are not on the
diagonal are polynomials with no degree zero term in
$e^{t_1},\ldots, e^{t_l}$,
\item[(h)] for any $1\le j\le m-2$, the block matrix $V_j$ is $(j+2)$-triangular and
its entries are polynomials with no degree zero term in
$e^{t_1},\ldots, e^{t_l}$.
\end{itemize}
By identifying the coefficients of powers of $\h$, the equation
(\ref{diffeq}) is equivalent to the system consisting of:
\begin{align}\label{sys1}d(V_0)V_0^{-1} &=
- \theta^{(1)}+[V_1, \omega]\end{align}
and
\begin{align}\label{sys2}
dV_1&=-\theta^{(2)}+[V_1,\theta^{(1)}]+[V_2,\omega]-V_1[V_1,\omega]\\
dV_i&=-\theta^{(i+1)}-\theta^{(i)}V_1-\ldots
-\theta^{(2)}V_{i+1}+[V_i,\theta^{(1)}]+[V_{i+1},\omega]-V_i[V_1,\omega]\nonumber
\end{align}
for $i\ge 2$.

It is  convenient to write a block matrix
$A=(A_{\alpha\beta})_{1\le \alpha, \beta\le m}$ as
$$A=A^{[-m]} +\ldots +A^{[-1]}+A^{[0]}+A^{[1]} +\ldots +A^{[m]}$$
where each block matrix $A^{[j]}$ is $j$-diagonal (i.e.
$A_{\alpha\beta}^{[j]}=0$ whenever $\beta-\alpha \ne j$). Then for
any two block matrices $A$ and $B$  we have:
$$(AB)^{[j]} = \sum_k A^{[k]}B^{[j-k]}, \quad [A,B]^{[j]}
=\sum_k[A^{[k]},B^{[j-k]}].$$

By (b), (c), (d) and (e) we can write:
\begin{align*}V_0=&I+V_0^{[2]}+V_0^{[3]}+\ldots + V_0^{[m]}\\
V_i=&V_i^{[i+2]} +V_i^{[i+3]}+\ldots + V_i^{[m]} \ \ (1\le i\le m-2)\\
\omega=&\omega^{[-1]} +\omega^{[0]} +\omega^{[1]}+\ldots
+\omega^{[m]}\\
\theta^{(i)}=&\theta^{(i),[i+1]}+\theta^{(i),[i+2]}+\ldots
+\theta^{(i),[m]} \ \ (1\le i\le m-1)
\end{align*}
In this way, the system (\ref{sys2}) is equivalent to:
\begin{align*} dV_1^{[j]} = -\theta^{(2),[j]}& +\sum_{3\le k\le
j-2}[V_1^{[k]}, \theta^{(1),[j-k]}]\\
 {} & +\sum_{4\le k\le j+1} [V_2^{[k]},\omega^{[j-k]}]\\
 {}& -\sum_{2\le k\le j-3}\sum_{3\le l\le
k+1}V_1^{[j-k]}[V_1^{[l]},\omega^{[k-l]}] \\
{}\\
 dV_i^{[j]} = -\theta^{(i+1),[j]} &-\sum_{3\le k\le
j-i-1}\theta^{(i),[j-k]}V_1^{[k]}\\
{} &-\sum_{i+3\le k\le j-3}\theta^{(2),[j-k]}V_{i+1}^{[k]}\\
{} & +\sum_{i+2\le k\le j-2}[V_i^{[k]},\theta^{(1),[j-k]}]\\
{} & +\sum_{i+3\le k\le j+1}[V_{i+1}^{[k]},\omega^{[j-k]}]\\
{} & -\sum_{2\le k\le j-i-2}\sum_{3\le l\le k+1} V_i^{[j-k]}
[V_1^{[l]},\omega^{[k-l]}],
\end{align*}
where $i\ge 2$. Define the total order on the matrices
$V_i^{[j]}$, $j\ge i+2$ as follows: $V_{i_1}^{[j_1]} <
V_{i_2}^{[j_2]}$ if and only if $j_1-i_1 <j_2-i_2$ or $j_1-i_1
=j_2-i_2$ and $j_1 <j_2$. We note that the system from above is of
the form:
$$dV_i^{[j]} = {\rm \ expression \ involving }\ V_{i'}^{[j']} >
V_i^{[j]},$$ for $i\ge 1$ and $j\ge i+2$. Because all coefficients
of the matrices $V_i^{[j]}$ are polynomials with no degree zero
term in $e^{t_1}, \ldots, e^{t_l}$, we deduce inductively ---
starting with $V_{m-2}^{[m]}$ ---  that there exists at most one
solution $V_i^{[j]}$, $i\ge 1$, $j\ge i+2$, of the system. It
remains to show that there exists at most one $V_0$ which
satisfies both the condition (f) and the equation (\ref{sys1}). If
$V_0'$ is another solution, then a simple calculation shows that
$$d(V^{-1}_0V_0' )=0.$$ By condition (f), the matrix
$V^{-1}_0V_0'$ has the diagonal $I$ and any  entry of it which is
not on the diagonal is a polynomial with no degree zero term in
$e^{t_1}, \ldots, e^{t_l}$. So $V^{-1}_0V_0'=I$, which means
$V_0=V_0'$.

\end{proof}

Now we can prove our main result:

{\it Proof of Theorem \ref{main}} Let $\star$ be a product with
the properties stated in Theorem \ref{main}. Consider the
isomorphism of $\D$-modules $$\phi:\D/\I\to H^*(G/B)\otimes
\bR[\{q_i\},\h]$$ given by Proposition \ref{isomorphism}. The
basis $\{[c_w]:w\in W\}$ of the right hand side induces the basis
$\{[\bar{C}_w]=\phi^{-1}([c_w]): w\in W\}$ of $\D/\I$ over
$\bR[\{Q_i\}, h]$. It is obvious that the latter satisfies the
hypotheses (i) and (iii) of Proposition \ref{unique}. We show that
it also satisfies (ii) by using the argument already employed in
the first part of the proof of Proposition \ref{unique}. Now from
Proposition \ref{unique}, we deduce that
$$[\bar{C}_w]=[\hat{C}_w],$$ for $w\in W$, where the basis
$\{[\hat{C}_w]:w\in W\}$ is induced by the actual quantum product
$\circ$.  Now, since $\phi$ is an isomorphism of $\D$-modules,
$\phi([\bar{C}_w])=[c_w]$ and $\phi(P_i)=[\lambda_i]$, we deduce
that the matrix of $[\lambda_i]\star$ with respect to the basis
$\{[c_w]:w\in W\}$ is the same as the matrix of $P_i$ with respect
to the basis $\{[\bar{C}_w]:w\in W\}$. Consequently we have
$$[\lambda_i]\star a=[\lambda_i] \circ a,$$
for all $a\in H^*(G/B) \otimes \bR[q_1,\ldots, q_l].$ Hence the
products $\star$ and $\circ$ are the same. \qed

\section{Quantization map for $Fl_n$}\label{last}

In the case $G=SL(n,\bC)$, the resulting flag manifold is $Fl_n$, which is 
the space of all complete flags in $\bC^n$. Borel's presentation (see eq. (\ref{borel})) in this case
reads
$$H^*(Fl_n) = \bR[\lambda_1,\ldots, \lambda_{n-1}]/(I_n)_{\ge 2},$$
where $(I_n)_{\ge 2}$ denotes the ideal generated by the nonconstant symmetric polynomials of degree at least 2
in the variables $$x_1:=-\lambda_1, x_2:=\lambda_1-\lambda_2,\ldots, x_{n-1}:=\lambda_{n-2}-\lambda_{n-1},
x_n:=\lambda_{n-1}.$$
Equivalently, we have    
$$H^*(Fl_n)=\bR[x_1,\ldots, x_n]/I_n$$
where $I_n$ denotes the ideal generated by the nonconstant symmetric polynomials of degree at least 1
in the variables $x_1,
\ldots, x_n.$ For any $k\in \{0,1,\ldots, n\}$ we consider the polynomials $e^k_0,\ldots, e^k_k$ in the variables
$x_1,\ldots, x_k$, which can be described by
$$\det \left[ \left(%
\begin{array}{ccccccc}
  x_1 &  0 &  \ldots & 0\\
  0 & x_2 & \ldots & 0 \\
  \ldots & \ldots & \ldots & \ldots \\
 0 & \ldots &  0 &  x_{k} \\
\end{array}%
\right) +\mu I_k \right] = \sum_{i=0}^n e_i^k\mu^{k-i} .
$$
For $i_1,\ldots, i_{n-1}\in \bZ$ such that $0\le i_j \le j$, we define
$$e_{i_1 \ldots i_{n-1}}=e_{i_1}^1\ldots e_{i_{n-1}}^{n-1}.$$
These are called the {\it standard elementary monomials.} 
It is known (see e.g. [Fo-Ge-Po, Proposition 3.4]) that the set
$\{[e_{i_1 \ldots i_{n-1}}] \ : \ 0\le i_j \le j\}$ is a basis of 
$H^*(Fl_n)$. 

We also consider the polynomials\footnote{These are the polynomials $E_i^k$ of
[Fo-Ge-Po].} $\hat{e}^k_0,\ldots, \hat{e}^k_k$ in the variables
$x_1,\ldots, x_k, q_1, \ldots, q_{k-1}$, which are described by
$$\det \left[ \left(%
\begin{array}{ccccccc}
  x_1 &  q_1 & 0 &  \ldots & 0\\
  -1 & x_2 &  q_2 &\ldots & 0 \\
  \ldots & \ldots & \ldots & \ldots & \ldots \\
  0 &\ldots & -1 & x_{k-1} & q_{k-1} \\ 
 0 & \ldots & 0 & -1 &  x_{k} \\
\end{array}%
\right) +\mu I_k \right] = \sum_{i=0}^k \hat{e}_i^k\mu^{k-i} .
$$
For $i_1,\ldots, i_{n-1}$ such that $0\le i_j \le j$, we define the 
{\it quantum standard elementary monomials}
$$\hat{e}_{i_1 \ldots i_{n-1}}=\hat{e}_{i_1}^1\ldots ,\hat{e}_{i_{n-1}}^{n-1}.$$
By a theorem of Ciocan-Fontanine  [Ci] (in fact Kim's theorem for $G=SL(n,\bC)$, see section 1), we have
the following isomorphism of $\bR[q_1,\ldots, q_{n-1}]$-algebras
\begin{equation}\label{quantizationmap}(H^*(Fl_n)\otimes\bR[q_1,\ldots, q_{n-1}],\circ)\simeq QH^*(Fl_n) := 
\bR[x_1,\ldots x_n,q_1,\ldots, q_{n-1}]/\langle \hat{e}_1^n,\ldots ,\hat{e}_{n}^{n}\rangle,
\end{equation}
which is canonical, in the sense that $[x_i]$ is mapped to $[x_i]_q$. According to [Fo-Ge-Po], we will
call this the {\it quantization map}.  
Since the conditions (\ref{deg}) and (\ref{hat}) are satisfied, we deduce that 
$\{[\hat{e}_{i_1 \ldots i_{n-1}}]_q \ : \ 0\le i_j \le j\}$ is a basis of 
$QH^*(Fl_n)$ over $\bR[q_1,\ldots, q_{n-1}]$. We also point out the obvious fact that 
$\{[{e}_{i_1 \ldots i_{n-1}}] \ : \ 0\le i_j \le j\}$ is a basis of 
$H^*(Fl_n) \otimes \bR[q_1,\ldots, q_{n-1}]$ over $\bR[q_1,\ldots, q_{n-1}]$. 
The goal of this section is to give a different proof to the following theorem of
Fomin, Gelfand, and Postnikov. 
\begin{theorem}\label{main3}{\rm (see [Fo-Ge-Po, Theorem 1.1])}. The quantization map described by equation
 (\ref{quantizationmap}) 
sends  $[e_{i_1\ldots i_{n-1}}]$ to $[\hat{e}_{i_1\ldots i_{n-1}}]_q.$
\end{theorem}

The main instrument of our proof is the $\D$-module $\D/\I$ defined in section 2. In this case
(i.e. $G=SL(n,\bC)$) we can describe it explicitly, as follows: 
$\D$ is  the (noncommutative) Heisenberg algebra  defined at the beginning of section 2, where $l=n-1$. 
The left ideal $\I$ of $\D$  is generated by $\E^n_1,\ldots, \E^n_{n-1}$, where
$$\det \left[ \left(%
\begin{array}{ccccccc}
  -P_1 &  Q_1 & 0 &  \ldots & 0\\
  -1 & P_1  -P_2  &  Q_2 &\ldots & 0 \\
  \ldots & \ldots & \ldots & \ldots & \ldots \\
  0 &\ldots & -1 & P_{n-2}  -P_{n-1} & Q_{n-1} \\ 
 0 & \ldots & 0 & -1 &  P_{n-1} \\
\end{array}%
\right) +\mu I_n \right] = \sum_{i=0}^n \E_i^n\mu^{n-i} .
$$
In fact we will need more general elements of $\D$, namely, for each $k\in \{1,\ldots, n-1\}$, we consider
the elements $\E_i^k$ of $\D$, with $0\le i\le k$, given by
 $$\det \left[ \left(%
\begin{array}{ccccccc}
  -P_1 &  Q_1 & 0 &  \ldots & 0\\
  -1 & P_1  -P_2  &  Q_2 &\ldots & 0 \\
  \ldots & \ldots & \ldots & \ldots & \ldots \\
    0 &\ldots & -1 & P_{k-2}  -P_{k-1} & Q_{k-1} \\ 
  0 &\ldots & 0 & -1 & P_{k-1} - P_k   \\ 
\end{array}%
\right) +\mu I_k \right] = \sum_{i=0}^k \E_i^k\mu^{k-i}.
$$
One can easily see that when we expand the determinant in the left hand side
of the last equation,
we will have no occurrence of  $P_jQ_j$ or $Q_jP_j$, $1\le j\le k-1$. This means that
the lack of commutativity of
$Q_j$ and $P_j$ creates 
no ambiguity in the definition of $\E^n_1,\ldots, \E^n_{n-1}$. We can also deduce that 
each of $\E^k_1,\ldots, \E^k_{k}$ is a linear combination of monomials in the variables $\{P_1,\ldots, P_{k},
 Q_1,\ldots, Q_{k-1}\}$,
with no ocurrence of $P_jQ_j$ or $Q_jP_j$ (i.e. the order of factors in each 
monomial is  not important). As a consequence, the following recurrence formula [Fo-Ge-Po, equation (3.5)] still
holds:
\begin{equation}\label{fgp}\E_{i}^k = \E_i^{k-1} + X_k \E_{i-1}^{k-1}+Q_{k-1}\E_{i-2}^{k-2},\end{equation}
where $X_k$ stands for $P_{k-1}-P_k$ and, by convention, $\E_j^k=0$, unless $0\le j \le k$.
It is worth mentioning the following commutation relations,
which will be used later:
 \begin{equation}\label{commutation}[X_k, \E_j^{l}]=0,\quad
[Q_k, \E_j^{l}]=0,\end{equation}
whenever $l\le k-1$.
We also  note that $\E_0^k=1$ and $\E_1^k = -P_k$ (where $P_n$ is by convention equal to 0).
We will prove the following result.
\begin{lemma} The elements $\E^k_1,\ldots, \E^k_{k-1}$ of $\D$ commute with each other.
\end{lemma}
\begin{proof} Consider the coordinates $s_0,\ldots, s_{k-1}$ on $\bR^k$. Following [Kim-Joe],
we consider the differential operators $D_j(\h\frac{\partial }{\partial s_0},
,\ldots,\h\frac{\partial }{\partial s_{k-1}}, e^{s_1-s_0}, \ldots, e^{s_{k-1}-s_{k-2}})$ given by  
 $$\det \left[ \left(%
\begin{array}{ccccccc}
  \h \frac{\partial }{\partial s_0} &  e^{s_1-s_0} & 0 &  \ldots & 0\\
  -1 & \h \frac{\partial }{\partial s_1} &  e^{s_2-s_1} &\ldots & 0 \\
  \ldots & \ldots & \ldots & \ldots & \ldots \\
    0 &\ldots & -1 & \h \frac{\partial }{\partial s_{k-2}} & e^{s_{k-1}-s_{k-2}} \\ 
  0 &\ldots & 0 & -1 & \h \frac{\partial }{\partial s_{k-1}}  \\ 
\end{array}%
\right) +\mu I_k \right] = \sum_{i=0}^k D_i^k\mu^{k-i}.
$$
By [Kim-Joe, Proposition 1], we have $[D_i^k, D_j^k]=0$ for all $0\le i,j\le k$. 
In order to prove our lemma, it is sufficient to note that if we make the change of coordinates
$$s_1-s_0=t_1, \ldots , s_{k-1}-s_{k-2} = t_{k-1}, -s_{k-1}=t_k,$$
we obtain
$$\h \frac{\partial }{\partial s_0} = - \h \frac{\partial }{\partial t_1}
=-P_1, \h \frac{\partial }{\partial s_1} =  \h \frac{\partial }{\partial t_1} - \h \frac{\partial }{\partial t_2}
=P_1-P_2, \ldots, \h \frac{\partial }{\partial s_{k-1}} =  
\h \frac{\partial }{\partial t_{k-1}} - \h \frac{\partial }{\partial t_{k}}
=P_{k-1}-P_{k},$$
where we have used the presentation of $\D$ given by $P_i=\h \frac{\partial }{\partial t_{i}}, Q_i=e^{t_i}$,
$1\le i \le n-1$.
\end{proof}
The following technical result will be needed later.
\begin{lemma}\label{commq}  We have 
\begin{equation}\label{firste} [\E_{j+1}^{k+1}, \E_i^k] = [\E_{i+1}^{k+1}, \E_j^k].\end{equation}
\end{lemma}
\begin{proof}
We  prove this by induction on $k \ge 0$. For $k=0$, the equation is obvious (by the convention
made above, we have 
$\E_0^j=0$). 
It  follows the induction step. We use the  recurrence formula (\ref{fgp}).
This gives
$$[ \E_{j+1}^{k+1}, \E_i^k] = [ \E_{j+1}^{k} + X_{k+1} \E_{j}^{k}+Q_{k} \E_{j-1}^{k-1}, \E_i^k]=
[Q_{k} \E_{j-1}^{k-1}, \E_i^k].$$ 
We continue by using again equation (\ref{fgp}) and obtain
 \begin{align*}{} &  [Q_{k} \E_{j-1}^{k-1}, \E_i^{k-1} + X_k \E_{i-1}^{k-1}+Q_{k-1}\E_{i-2}^{k-2}]
 \\{} &  =  [Q_k, X_k] \E_{i-1}^{k-1}\E_{j-1}^{k-1} + [Q_k\E_{j-1}^{k-1}, Q_{k-1} \E_{i-2}^{k-2}]
  \\{} &  =  [Q_k, X_k] \E_{i-1}^{k-1}\E_{j-1}^{k-1} + Q_k[\E_{j-1}^{k-1}, Q_{k-1} \E_{i-2}^{k-2}]
 \\ {} & = [Q_k, X_k] \E_{i-1}^{k-1}\E_{j-1}^{k-1} + Q_k[\E_{j-1}^{k-1},\E_{i}^k - \E_i^{k-1} - X_k \E_{i-1}^{k-1}]  
\\{} &  =[Q_k, X_k] \E_{i-1}^{k-1}\E_{j-1}^{k-1} +Q_k([\E_{j-1}^{k-1}, \E_{i}^k] - 
[\E_{j-1}^{k-1}, X_k \E_{i-1}^{k-1}])\\ {} & =  [Q_k, X_k] \E_{i-1}^{k-1}\E_{j-1}^{k-1} +Q_k[\E_{j-1}^{k-1}, \E_{i}^k]
\end{align*}
Here we have used the commutation relations (\ref{commutation}) several times.   
 Similarly, we obtain
 $$[\E_{i+1}^{k+1}, \E_j^k] = 
[Q_k, X_k] \E_{j-1}^{k-1}\E_{i-1}^{k-1} +Q_k [\E_{i-1}^{k-1}, \E_{j}^{k}].$$
We use the induction hypothesis  to finish the proof.
\end{proof}

Now we are ready to prove Theorem \ref{main3}.

\noindent {\it Proof of Theorem \ref{main3}.} Let $\omega_k$ denote the matrix of multiplication by $[y_k]_q$ with respect to
 the basis $\{[\hat{e}_{i_1 \ldots i_{n-1}}]_q \ : \ 0\le i_j \le j\}$  of 
$QH^*(Fl_n)$ (see equation (\ref{quantizationmap})). More precisely, the entries of $\omega_i$ are polynomials in 
$q_1,\ldots, q_{n-1}$, determined by
\begin{equation}\label{coefficients}  [y_k]_q [ \hat{e}_{i_1 \ldots i_{n-1}}]_q =
\sum_{l_1,\ldots, l_{n-1}} \omega_k^{i_1 \ldots i_{n-1}, l_1 \ldots l_{n-1}}[\hat{e}_{l_1 \ldots l_{n-1}}]_q.
\end{equation}
According to Corollary \ref{proper}, it is sufficient to show that
\begin{equation}\label{partiali}\frac{\partial}{\partial t_i}\omega_j =\frac{\partial}{\partial t_j}\omega_i,
\end{equation}
for $1\le i,j\le n-1$, where as usually, we use the convention $q_i=e^{t_i}$. 
For $i_1,\ldots, i_{n-1}$ such that $0\le i_j \le j$, we consider
$$\E_{i_1 \ldots i_{n-1}}:=\E_{i_1}^1 \E_{i_2}^2 \ldots \E_{i_{n-1}}^{n-1}.$$
In order to prove equation (\ref{partiali}), it is sufficient to
prove the following claim.

\noindent {\it Claim.}  In $\D/\I$ 
we have 
\begin{equation}\label{claim3}[P_k] [ \E_{i_1 \ldots i_{n-1}}] =\sum_{l_1,\ldots, l_{n-1}} \Omega_k^{i_1 \ldots i_{n-1}, 
l_1 \ldots l_{n-1}}[\E_{l_1 \ldots l_{n-1}}],\end{equation}
where each $\Omega_k^{i_1 \ldots i_{n-1}, l_1 \ldots l_{n-1}}$ 
is obtained from $\omega_k^{i_1 \ldots i_{n-1}, l_1 \ldots l_{n-1}}$
 by  the modification $Q_i\mapsto q_i$. 
 
 Indeed, if we make the usual identifications $P_k=\h\frac{\partial}{\partial t_k}$, $Q_k=e^{t_k}$,
$1\le k\le n-1$, then (\ref{claim3}) implies that the connection
 $$d+ \sum_{k=1}^{n-1}\frac{1}{\h}\Omega_kdt_k$$ is flat (see e.g. [Gu, Proposition 1.1]) for all 
values of $\h$, which implies (\ref{partiali}).
The proof of the claim relies on  a noncommutative version of the quantum straightening algorithm
of Fomin, Gelfand, and Postnikov [Fo-Ge-Po]. The key equation is the following.
\begin{equation}\label{straightening}\E_i^k\E_{j+1}^{k+1} + \E_{i+1}^k \E_j^k +Q_k\E_{i-1}^{k-1}\E_j^k = \E_j^k \E^{k+1}_{i+1} 
+\E_{j+1}^k\E_i^k +Q_k\E_{j-1}^{k-1}\E_i^k.\end{equation}
We note that this is the same as   equation (3.6) in [Fo-Ge-Po]. The difference is that here 
we work in the algebra $\D$, which is not commutative, so it is not {\it a priori} clear that
(\ref{straightening}) still holds.  In order to prove it,
we use equation (\ref{fgp}) twice and obtain:
$$(\E_{j+1}^{k+1}-\E_{j+1}^k)\E_i^k = (X_{k+1}\E_j^k +Q_k\E_{j-1}^{k-1})\E_i^k,$$
and
$$ (\E_{i+1}^{k+1}-\E_{i+1}^k)\E_j^k = (X_{k+1}\E_i^k +Q_k\E_{i-1}^{k-1})\E_j^k.$$ 
If we subtract the second equation from the first one, we obtain:
$$\E_{i+1}^{k+1}\E_j^k - \E_{j+1}^{k+1}\E_i^k=  \E_{i+1}^k\E_j^k - \E_{j+1}^k\E_i^k 
+Q_k(\E_{i-1}^{k-1}\E_j^k  - \E_{j-1}^{k-1}\E_i^k).$$
Now the left hand side can be written as
$$\E_j^k\E_{i+1}^{k+1} - \E_i^k\E_{j+1}^{k+1} + [\E_{i+1}^{k+1},\E_j^k] - [\E_{j+1}^{k+1},\E_i^k]
= \E_j^k\E_{i+1}^{k+1} - \E_i^k\E_{j+1}^{k+1} ,$$
where we have used Lemma \ref{commq}. Equation (\ref{straightening}) has been proved.
Now we can use it exactly like in the commutative situation, described
in [Fo-Ge-Po], in order to obtain the expansion of the product of 
$P_k= -\E_1^k$ and $\E_{i_1 \ldots i_{n-1}} = \E_{i_1}^1 \ldots  \E_{i_{n-1}}^{n-1}$.
More precisely, we begin with 
$$P_k \E_{i_1 \ldots i_{n-1}} = 
\E_{i_1}^1 \ldots \E_{i_{k-1}}^{k-1} P_k \E_{i_k}^k \E_{i_{k+1}}^{k+1} \ldots  \E_{i_{n-1}}^{n-1} 
= -\E_{i_1}^1 \ldots \E_{i_{k-1}}^{k-1} \E_1^k \E_{i_k}^k \E_{i_{k+1}}^{k+1} \ldots  \E_{i_{n-1}}^{n-1} ,$$
and then we use (\ref{straightening}) repeatedly. The resulting coefficients in the final expansion
will be the same as in the commutative situation. This finishes the proof of the claim, and also of 
Theorem \ref{main3}.  \hfill $\square$      

\bibliographystyle{abbrv}

\begin{thebibliography}{Fo-Ge-Po}

\bibitem[Am-Gu]{Am-Gu} A. Amarzaya and M. A. Guest, {\em Gromov-Witten invariants of flag
manifolds, via D-modules}, Jour. London Math. Soc. (2) {\bf 72} (2005), 121--136

\bibitem[Be-Ge-Ge]{Be-Ge-Ge} I. N. Bernstein, I. M. Gelfand and S. I. Gelfand,
{\it Schubert cells and cohomology of the space $G/P$}, Russian Math. Surveys {\bf  28} (1973),
1--26

\bibitem[Bo]{Bo} A. Borel, {\it Sur la cohomologie des espaces fibr\'es
principaux et des espaces homog\`enes des groupes de Lie
compacts},  Ann. of Math.  {\bf 57} No. 2 (1953), 115--207

\bibitem[Ci]{Ci} I. Ciocan-Fontanine,
{\it The quantum cohomology ring of flag varieties}, Trans. Amer. Math. Soc. {\bf 351} (1999), no. 7, 2695--2729

\bibitem[Du]{Du} B. Dubrovin, {\it The geometry of 2D topological field theories},
 Integrable Systems and Quantum Groups, Lecture Notes in
 Mathematics, Vol. 1620, Springer-Verlag, New York, 1996, 120--348

\bibitem[Fo-Ge-Po]{Fo-Ge-Po}  S. Fomin, S. Gelfand, and A. Postnikov, {\em Quantum
Schubert polynomials},  J. Amer. Math. Soc.,  {\bf 10} (1997),
565--596


\bibitem[Fu-Pa]{Fu-Pa}  W. Fulton and R. Pandharipande,  {\em Notes on stable
maps and quantum cohomology}  Algebraic geometry---Santa Cruz
1995,  Proc. Sympos. Pure Math., 62, Part 2, editors J. Kollar, R.
Lazarsfeld and D.R. Morrison,  1997,  45--96

\bibitem[Fu-Wo]{Fu-Wo} W. Fulton and C. Woodward, {\it On the
quantum product of Schubert classes}, J. Algebraic Geom. {\bf 13} (2004), 641--661

\bibitem[Gu]{Gu}  M. A. Guest, {\em Quantum cohomology via D-modules}, Topology {\bf 44}
(2005), 263--281


\bibitem[Ir]{Ir} H. Iritani, {\em Quantum D-module and equivariant
Floer theory for free loop spaces}, preprint {\tt math.DG/0410487}

\bibitem[Kim]{Kim} B. Kim, {\em Quantum cohomology of flag manifolds
$G/B$ and quantum Toda lattices},  Ann. of Math. {\bf 149}
(1999), 129--148

\bibitem[Kim-Joe]{Kim-Joe} B. Kim and D. Joe, 
{\it Equivariant mirrors and the Virasoro conjecture for flag manifolds}, 
Int. Math. Res. Not. {\bf 15} (2003),  859--882


\bibitem[Ma1]{Ma1} A.-L. Mare, {\em  On the theorem of Kim
concerning $QH^*(G/B)$}, Integrable systems, topology and physics,
editors M. Guest, R. Miyaoka and Y. Ohnita, Contemp. Math. 309, Amer.
Math. Soc. (2002), 151-163

\bibitem[Ma2]{Ma2} A.-L. Mare, {\em Polynomial representatives of Schubert classes
in $QH^*(G/B)$}, Math. Res. Lett. {\bf 9} (2002), 757--770

\bibitem[Ma3]{Ma3} A.-L. Mare, {\em
Relations in the quantum cohomology ring of $G/B$}, Math. Res.
Lett. {\bf 11} (2004), 35--48

\bibitem[Ma4]{Ma4} A.-L. Mare, {\em The combinatorial quantum cohomology ring of $G/B$}, 
Jour. Alg. Comb. {\bf 21} (2005), 331--349

\bibitem[Pe]{Pe} D. Peterson, {\em Lectures on quantum cohomology
of $G/P$}, M.I.T. 1996

\bibitem[Si-Ti]{Si-Ti} B. Siebert and G. Tian, {\it On quantum cohomology
rings of Fano manifolds and a formula of Vafa and Intriligator},
Asian J. Math. {\bf 1} (1997), 679-695

\end{thebibliography}

\end{document}